\documentclass{article}
\usepackage{amssymb}

\begin{document}

\newcommand{\C}{\mathbb{C}}
\newcommand{\R}{\mathbb{R}}
\newcommand{\N}{\mathbb{N}}
\newcommand{\K}{\mathbb{K}}
\newcommand{\Q}{\mathbb{Q}}
\newcommand{\pr}{\mbox{\rm pr}}
\newcommand{\af}{\mbox{\rm Aff}(\mathbb{R}^n)}

\newtheorem{theorem}{Theorem}
\newtheorem{proposition}[theorem]{Proposition}
\newtheorem{defi}[theorem]{Definition}
\newtheorem{corollary}[theorem]{Corollary}
\newtheorem{lemma}[theorem]{Lemma}
\newtheorem{remark}[theorem]{Remark}
\newcommand{\dem}{\noindent{\bf Proof. }}
\newcommand{\qed}{\hfill $\square$}

\title{Invariant quadrics and orbits for a family of rational systems of difference equations} 

\author{Ignacio Bajo}

\maketitle

\begin{abstract}
We study the existence of invariant quadrics for a class of  systems of difference equations in ${\mathbb R}^n$ defined by linear fractionals sharing denominator. Such systems can be described in terms of some square matrix $A$ and we prove that there is a correspondence between non-degenerate invariant quadrics and solutions to a certain matrix equation involving $A$. We show that if  $A$ is semisimple and the corresponding system admits non-degenerate quadrics, then every orbit of the dynamical system is contained either in an invariant affine variety or in an invariant quadric.  
\end{abstract}

\section{Introduction}  

The determination of geometric invariants of a system of difference equations is, in general, a difficult task. In this paper we will find a family of quadratic varieties which remain invariant for a class of systems of difference equations in $\R^n$, $n\ge 2$, of the type
\begin{equation}\label{F}
 X(k+1)=(F_1(X(k)),\dots,F_n(X(k))),\quad X(k)\in\R^n\end{equation}
 where the maps $F_i$ are linear fractionals sharing denominator:
$$F_i(X)=\frac{a_{i1}x_1+a_{i2}x_2+\cdots+a_{in}x_n+c_{i}}{b_{1}x_1+b_{2}x_2+\cdots+b_{n}x_n+d},\, \ i=1,2,\dots , n,$$
where $X=(x_1,x_2,\cdots,x_n)$ and all the involved parameters are real.
Such kind of rational systems have been treated in \cite{b-l} where global periodicity properties were  studied. Further, AlSharawi and Rhouma \cite{Al-R} studied a biological model given by systems of rational difference equations with a common denominator. The key fact for the study of those systems is that they can be written in certain matricial form by the use of homogeneous coordinates. Explicitly, if one denotes by $q$ the mapping given by $q(a_1,a_2,\dots,a_{n+1})=(a_1/a_{n+1},a_2/a_{n+1},\dots,a_n/a_{n+1})$ for $(a_1,a_2,\dots ,a_{n+1})\in {\mathbb R}^{n+1}$ with $a_{n+1}\neq 0$ and $\ell: {\mathbb R}^n\to {\mathbb R}^{n+1}$ is the transformation to homogeneous coordinates  given by $\ell ( a_1,\dots,a_n)=(a_1,\dots,a_n,1)^T$, where $M^T$ stands for the transposed of a matrix $M$, then the system can be written in the  form $X(k+1)=q\circ A\circ\ell (X(k))$  
for the square matrix 
of order $(n+1)$ given  by
\begin{equation}\label{A}
A=\left(\begin{array}{c|c} A_1 &C^T\\ \hline B & d\end{array}\right),
\end{equation}
where $A_1=(a_{ij})$, $B=(b_1,\dots , b_n)$ and $C=(c_1,\dots , c_n).$

The following result of \cite{b-l} will be frequently used:

\begin{lemma}
With the notations above, we have:
\begin{enumerate}
\item[{\rm (a)}] $q(\ell(x))=x$ for all $x\in\R^n.$
\item[{\rm (b)}] $\ell(q(a))=(1/a_{n+1})a$ for all $a=(a_1,\dots, a_{n+1})$ such that $a_{n+1}\ne 0$.
\item[{\rm (c)}] If $a,b\in\R^{n+1}$ are such that $q(a)$ and $q(b)$ exist, then $q(a)=q(b)$ if and only if $a=\lambda b$ for some $\lambda\in\R,\, \lambda\neq 0$.
\item[{\rm (d)}] If $A$ is an $(n+1)\times (n+1)$ matrix, then $q(A\ell(q(a)))=q(Aa)$ whenever both members exist.
\end{enumerate}
\end{lemma}

Since we will be interested in orbits and invariant sets, we recall the following: 

Let $F:\Omega\subset\R^n\to\R^n$ be a map defined on certain non-empty set $\Omega$ and $X(0)\in\R^n$. In we denote by $F^k$ the $k$-th power under composition of $F$, then the {\it orbit} of $X_0$ is the sequence $\{F^k(X_0)\}_{k\in\N}$ or, equivalently, the solution of the system $X(k+1)=F(X(k))$ with initial condition $X(0)=X_0$. If for some $k\in\N$ the corresponding power $F^k(X_0)$  is not defined we say that the orbit is {\it uncomplete}, and that it is {\it complete} or a {\it full orbit} otherwise. 

A subset $\mathcal{S}\subset\R^n$ is said to be {\it $F$-invariant} if $F(\Omega\cap\mathcal{S})\subset\mathcal{S}$. This means that every orbit starting at a point of $\mathcal{S}$ remains in $\mathcal{S}$.

For our system (\ref{F}), if we denote  $F=(F_1,\dots,F_n)$, it is obvious that the map $F$ is defined on $\R^n\setminus \mathcal{PF}$, where
$$\mathcal{PF}=\{\ (x_1,x_2,\cdots,x_n)\in\R^n\ :\ b_{1}x_1+b_{2}x_2+\cdots+b_{n}x_n+d=0\ \}$$
is the, so called, {\it principal forbidden hyperplane}. Note that if $A$ is the matrix such that $F=q\circ A\circ\ell$ then we actually have
$$\mathcal{PF}=\{\ X\in\R^n\ :\ (0,\dots,0,1)A\ell(X)^T=0\ \}.$$
Further, one may easily see \cite{b-l} that $F^k=q\circ A^k\circ\ell$ and, hence, the orbit starting at a point $X(0)$  is complete if and only if $\pr (A^k\circ\ell(X(0)))\ne 0$ for all $k\in\N$, where $\pr:\R^{n+1}\to\R$ is the projection on the $(n+1)$ component
$\pr(x_1,x_2,\cdots,x_{n+1})=x_{n+1}$.
From now on, we will denote $U_0$ the inverse image of $\{0\}$ by $\pr$, this is to say
$$U_0=\{ u\in\R^{n+1}\ :\ \pr(u)=0\ \}.$$
 
\medskip

In the sequel, we will use the term {\it quadric in } $\R^{n}$ for the locus in $\R^{n}$ of zeros of a quadratic polynomial. We will write the quadratic equation in homogeneus coordinates and, therefore, a quadric in $\R^{n}$ will be given by the set
$$\mathcal{Q}(M)=\{X\in\R^n\ :\ \ell(X)^TM\ell(X)=0\}$$
for certain symmetric matrix $M=M^T\in \R^{(n+1)\times (n+1)}$. If $M$ is non-singular, we say that the quadric is {\it non-degenerate}. As usual, quadrics in $\R^{2}$ will be called {\it conics}.

It was seen in \cite{bfp} for some particular cases in $\R^2$ that each orbit  is contained in a certain conic when the corresponding matrix $A$ is semisimple (i.e. diagonalizable as a complex matrix) and all its eigenvalues have the same modulus. This suggested the idea that a similar result could be also verified for higher dimensional cases and, moreover, that the existence of invariant algebraic varieties of degree two could be proved. The aim of this paper is to study, in the general $n$-dimensional case, for which matrices does this occurs and to derive some geometrical considerations regarding the orbits of the associated system. In particular, we prove that when the matrix is semisimple and similar to its inverse multiplied by $\pm 1$, the existence of non-degenerate quadrics can be guaranteed and every orbit lies either in an invariant affine variety or in an invariant quadric. 

From now on, we will denote the spectrum of a real matrix $A$ by $\mbox{Sp}(A)$. If  $M_1,M_2,\dots,M_k$ are square matrices, we will denote the corresponding block diagonal matrix by $\mbox{diag}(M_1,M_2,\dots,M_k).$

\section{Existence of invariant quadrics}

 We start with the following result on  the existence of certain $F$-invariant varieties for our rational system (\ref{F}). Notice that if $U$ is a subspace of $\R^{n+1}$ not contained in $U_0$, then one can choose a basis $\{u_1,\dots,u_k\}$ of $U$ such that $\pr(u_1)=1$ and $\pr(u_i)=0$ for $i\ge 2$. It is then clear that  $u=\sum_{i=1}^k\xi_iu_i\in U\setminus (U\cap U_0)$ if and only if $\xi_1\ne 0$ and that $\pr(u)=\xi_1$ and, therefore,  $q(u)=q(u_1)+\sum_{i=2}^k(\xi_i/\xi_1)X_i$, where $X_i=q(u_i+\ell(0))$ is the projection of $u_i$ on the first $n$ components. Since $\{X_i\}_{i=2}^k$ is a linearly independent set, we immediately get that   the set $\mathcal{S}_U=q(U\setminus (U\cap U_0))$ is an affine variety of dimension $(k-1)$.

\begin{proposition} \label{invpl} Let $A\in\R^{(n+1)\times(n+1)}$ be an arbitrary matrix and $F=q\circ A\circ\ell$. Let $U$ be a  subspace of $\R^{n+1}$  not contained in $U_0$ and suppose that the affine variety $\mathcal{S}_U=q(U\setminus (U\cap U_0))$ is not cointained in $\mathcal{PF}$.

The subspace $U$ is $A$-invariant if and only if  $\mathcal{S}_U$ is $F$-invariant. 
\end{proposition}
\dem Let us suppose first that $U$ is $A$-invariant and consider $X\in\mathcal{S}_U$, $X\not\in\mathcal{PF}$. We then have that $X=q(u)$ for some $u\in U$ such that $\pr(u)\ne 0$ and that $A(\ell(q(u)))\not\in U_0$. Since $\ell(q(u))=(1/\pr(u))u$, this means that $Au\not\in U_0$. Thus, we have
$$F(X)=q((1/\pr(u)Au)=q(Au)$$
where $Au\in U$ but $Au\not\in U_0$ and hence $F(X)\in\mathcal{S}_U$.

In order to prove the converse, notice that, since $\mathcal{S}_U$ is $F$-invariant, when $X\in\mathcal{S}_U$ and $F(X)$ exists, one has that there exists $\tilde u\in U$ such that $q(A(\ell(X)))=q(\tilde u)$, which implies that $A\ell(X)=\lambda \tilde u\in U$ for some $\lambda\in\R$. Recalling that $\mathcal{S}_U$ is not contained in $\mathcal{PF}$, we can always find an element $X_1\in \mathcal{S}_U$ verifying such conditions. Thus, if $u_1=\ell (X_1)$ then $u_1,Au_1\in U$. Now, it is easy to prove that for every $u\in U$ one can find $\alpha_1,\alpha_2\in\R$, $\alpha_2\ne 0$, such that $\pr(\alpha_1u_1+\alpha_2u)\neq 0$ and $\pr(\alpha_1Au_1+\alpha_2Au)\neq 0$, this meaning that $X=q(\alpha_1u_1+\alpha_2u)$ is well defined and $F(X)$ exists. Therefore, $\ell(q (\alpha_1Au_1+\alpha_2Au))\in U$, from where it is straightforward to prove that $Au\in U$.\qed

\begin{remark} {\em  In particular, if $A$ admits a real eigenvector outside $U_0$, then the $q$-projection of its linear span provides a fixed point of the rational system. It should be noticed that, since $n\ge 2$, when  the matrix $A$ is semisimple, it always admits, at least, one proper $A$-invariant subspace not contained in $U_0$ and, thus, there always exists a proper $F$-invariant affine variety.}\end{remark}

Although the proposition above gives a nice  result of existence of $F$-invariant varieties, it is obvious that not every orbit of the map $F$ lies completely inside one of those. Our aim is to find in which cases for each $X_0\not\in\mathcal{PF}$ one can always find an invariant quadric containing the whole orbit of $X_0$. 

When the matrix $A$ is singular, the $n$-dimensional rational system defined by $F=q\circ A\circ\ell$ can be reduced to a rational system of the same type in $(n-1)$ dimensions. Therefore from now on we will always consider that the matrix $A$ is invertible. Moreover, since for every matrix $A$ and $0\ne\lambda\in\R$ one has that $q\circ A\circ\ell=q\circ \lambda A\circ\ell$, we can  consider always that $\mbox{det}(A)=\pm 1.$ 

\begin{proposition} Let $F=q\circ A\circ\ell$ where $A\in\R^{(n+1)\times(n+1)}$ and let $M$ be an indefinite symmetric matrix such that $\mbox{\rm rank} (M)>2$. 

The quadric $\mathcal{Q}(M)$ is $F$-invariant if and only if $A^TMA=\mu M$ for some non zero $\mu\in\R$. 
\end{proposition}
\dem  Suppose that $X\in\R^n$ verifies that $\ell(X)^TM\ell(X)=0$ and  $\pr(A\ell(X))\ne 0$. If we put $\gamma=1/\pr(A\ell(X))$ then we have
$\ell(F(X))=\ell(qA\ell(X))=\gamma A\ell(X)$ and, hence,
$\ell(F(X))^TM\ell(F(X))=\gamma^2\ell(X)^TA^TMA\ell(X).$
This clearly shows  that if $A^TMA=\mu M$ for a certain $\mu\ne 0$, then $\mathcal{Q}(M)$ is $F$-invariant. 

For the converse, we fist remind  that if $M$ and $S$ are indefinite symmetric matrices with rank larger than 2 such that $\mathcal{Q}(M)=\mathcal{Q}(S)$ then they are proportional \cite{vinberg}. Thus, if we see that $\mathcal{Q}(M)=\mathcal{Q}(S)$ for $S=A^TMA$ then the equality $A^TMA=\mu M$ will follow immediatly. The main problem in order to see that $\mathcal{Q}(M)=\mathcal{Q}(S)$ is that, when $\mathcal{Q}(M)$ is $F$-invariant we only have, in principle, the equivalence between $\ell(X)^TM\ell(X)=0$ and $\ell(X)^TA^TMA\ell(X)=0$ for those vectors such that $\pr(A\ell(X))\ne 0$. Let us see that actually, the equivalence also holds when $\pr(A\ell(X))= 0$. Denote by $G:\mathcal{Q}(M)\to\R$ the continuous map $G(X)=\ell(X)^TA^TMA\ell(X)$. The set of points $X\in\mathcal{Q}(M)$ such that $\ell(X)^TA^TMA\ell(X)\ne 0$ is actually $G^{-1}(\R\setminus\{0\})$ and, hence, it is an open subset of $\mathcal{Q}(M)$. Since $M$ is indefinite, we can assure that $\mathcal{Q}(M)$ is a manifold of dimension $n-1$ and, thus, so is the open subset $G^{-1}(\R\setminus\{0\})$ if it is non-empty. Further, since $\ell(X)^TA^TMA\ell(X)\ne 0$ can only occur whenever $\pr(A\ell(X))=0$, we have that $G^{-1}(\R\setminus\{0\})$ is contained in the hyperplane $\pr(A\ell(X))=0$. Thus, $\mathcal{Q}(M)$ contains an open set contained in a hyperplane and, as a consequence, it must contain the whole hyperplane. But this is not possible since $\mbox{\rm rank} (M)>2$. This proves that $G^{-1}(\R\setminus\{0\})$ must be the empty set and $\mathcal{Q}(M)=\mathcal{Q}(S)$.\qed

\begin{remark}{\em First notice that the condition on the indefiniteness of $M$ and that on its rank ensure that $\mathcal{Q}(M)$ is not an affine variety and does not contain a hyperplane, which are essential in the proof. 

When $A^TMA=\mu M$ holds for a symmetric matrix $M$ and $\mu<0$ then the signatures of $M$ and $\mu M$ are equal and, hence,  must be of the form $(p,p)$ for some $p\le (n+1)/2$.

It is well known \cite{hj2,laub} that if $\mu\in\R$ is an eigenvalue of the Kronecker product $A^T\otimes A^T$, then there exist non-null solutions to the matrix equation $A^TYA=\mu Y$ and, therefore, if $Y$ is not skew-symmetric, one can always find a symmetric matrix $M=Y+Y^T$ verifying $A^TMA=\mu M$. Such matrix $M$ is not, in general, invertible. However,
from a dynamical point of view, is seems more interesting to find non-degenerate quadrics; for instance, the existence of invariant ellipsoids guarantees  the boundedness of the orbits they contain. Thus we will focus our attention in the case of an invertible symmetric matrix $M$. }\end{remark}

\begin{corollary} Let $F=q\circ A\circ\ell$ where $A\in\R^{(n+1)\times(n+1)}$ with $\det(A)=\pm 1$ and let $M$ be an invertible indefinite  symmetric matrix. 

The quadric $\mathcal{Q}(M)$ is $F$-invariant if and only if $A^TMA=\varepsilon M$ for some $\varepsilon\in\{1,-1\}$. As a consequence, the matrix $A$ must be similar to $\varepsilon A^{-1}$.
\end{corollary}
\dem Clearly, $A^TMA=\varepsilon M$ implies that $\det(A)^2\det(M)=\varepsilon^{n+1}\det(M)$ and thus $\varepsilon\in\{1,-1\}$ because $\det(A)=\pm 1$ and $\det(M)\ne 0$.

The second part of the statement is almost trivial since, when $M$ is invertible, the condition  obviously becomes $ A^T=\varepsilon MA^{-1}M^{-1}$, which imposes that $A^T$ (and, hence, $ A$) is similar to $\varepsilon A^{-1}$.\qed

\begin{remark}{\em Notice that the case $\varepsilon =-1$ is only possible if $n+1$ is even and, in such case, the signature of $M$ must be $((n+1)/2,(n+1)/2)$. 

The fact that $A$ is similar to $\varepsilon A^{-1}$ means that if $\lambda\in\C$ is an eigenvalue of $A$, so is $\varepsilon\lambda^{-1}$ and the Jordan block decompositon for both eigenvalues is the same. 
Unfortunately, although that condition guarantees the existence of an invertible matrix $Y$ such that $A^TYA=\varepsilon Y$, it does not assure that $Y$ can be taken symmetric. For example, if one considers the matrix
$$A=\left(\begin{array}{ccc} 1& 0& 0\\1& 1 &0\\0 &0 & 1\end{array}\right)\, ,$$
then the only symmetric matrices $M$ verifying $A^TMA=M$ are of the form
$$M=\left(\begin{array}{ccc} \alpha & 0& \beta\\0& 0 &0\\ \beta &0 & \gamma\end{array}\right).$$
In the case $\varepsilon =-1$ one can also find examples in which invertible solutions to  $A^TYA=- Y$ exist but none of them is symmetric. This is the case, for instance, of the matrix
$$A=\left(\begin{array}{rrrr} 0&1& 0& 0\\-1& 0 & 0&0\\1&0 &0 & 1\\0&1 &-1 & 0\end{array}\right)\, .$$

Thus, it will be important to determine further conditions to ensure the existence of invertible symmetric solutions to the matrix equation. Besides, one may be interested in determining in some way how large is the set of invertible symmetric matrices $M$ which verify $A^TMA=\varepsilon M$ for a given matrix $A$. A reasonable  measure for this seems the maximum number of independent matrices verifying such conditions. In Theorem \ref{existence} below we will see that if $A$ is semisimple and similar to $\varepsilon A^{-1}$, then such number is not null and we will determine it explicitely. We first prove  some previous results.
}\end{remark}

\begin{proposition} \label{red-co} Consider $\varepsilon=\pm 1$ and for $A\in\R^{(n+1)\times (n+1)}$ let us  put
 $$\mathcal{C}_\varepsilon (A)={\mathbb R}\mbox{\rm -span}\{ M\in \R^{(n+1)\times (n+1)}\, :\,\det(M)\ne 0,\, M=M^T,\ A^TMA=\varepsilon M\}.$$ 
 \begin{enumerate}
 \item[(1)] If $\mathcal{C}_\varepsilon (A)\ne\{0\}$, then
 $$\mathcal{C}_\varepsilon (A)=\{ M\in \R^{(n+1)\times (n+1)}\, :\, M=M^T,\ A^TMA=\varepsilon M\}.$$
 \item[(2)] The set $\mathcal{C}_\varepsilon (A)$ is non-null if and only if there exists an invertible matrix $N_0\in \C^{(n+1)\times (n+1)}$ such that  $N_0=N_0^T$ and $A^TN_0A=\varepsilon N_0$ and, in such case,
 $$\mbox{\rm dim}\,\mathcal{C}_\varepsilon (A)=\mbox{\rm dim}_\C\,\{N\in \C^{(n+1)\times (n+1)}\, :\, N=N^T,\ A^TNA=\varepsilon N\}$$
 where $\mbox{\rm dim}_\C(V)$ denotes the complex dimension of a complex  space $V$.
 \end{enumerate}
\end{proposition}
\dem To see the first assertion, it suffices to prove that if there exists $M_0\in \mathcal{C}_\varepsilon (A)$ then every singular matrix $M$ such that $M=M^T$ and $A^TMA=\varepsilon M$ is actually in $\mathcal{C}(A)$. Choose such a matrix $M$ and consider the polynomial in  $\mu\in\R$, given by $\det (M_0+ \mu M)\ne 0$, which is not identically zero since its value in $\mu=0$ is $\det (M_0)\ne 0$. Since  a non-null polynomial has a finite number of roots, we may find $\mu_0\in\R$, $\mu_0\ne 0$ such that $\det(M_0+\mu_0 M)\ne 0.$ Thus,  $M=\mu_0^{-1}( M_0+ \mu_0 M)-\mu_0^{-1} M_0\in \mathcal{C}(A)$.

In order to prove (2), note that if $N_0$ is  complex symmetric and verifies $A^TN_0A=\varepsilon N_0$ then its real and imaginary parts $M_1,M_2$ are real symmetric and verify $A^TM_iA=\varepsilon M_i$ for $i=1,2$. As before, if $N_0$ is invertible, then the polynomial in $\mu$ $\det(M_1+\mu M_2)$ is non zero since its value in the imaginary unit is $\det (N_0)\ne 0$ and we may find $\mu_0\in\R$ such that $\det(M_1+\mu_0 M_2)\ne 0.$ Thus, the matrix $M=M_1+\mu_0 M_2$ is real symmetric and invertible and clearly verifies $A^TMA=\varepsilon M$. The equality for the dimensions follows at once from (1) and the fact that
$$\{N\in \C^{(n+1)\times (n+1)}\, :\, N=N^T,\ A^TNA=\varepsilon N\}$$
is just the complexification of the set of real symmetric matrices $M$ verifying $A^TMA=\varepsilon M$.\qed

\begin{lemma} Let us choose $\varepsilon\in\{1,-1\}$ and  consider $A=\mbox{\rm diag} (A_1,A_2)$ where $A_1,A_2$ are invertible complex matrices of arbitrary orders  such that $\mbox{\rm Sp}(A_1)\cap\mbox{\rm Sp}(\varepsilon A_2^{-1})=\emptyset$.

A complex matrix $R$ verifies $A^TRA=\varepsilon R$ if and only if there exist $R_1,R_2$ such that $R=\mbox{\rm diag} (R_1,R_2)$ and $A_i^TR_iA_i=\varepsilon R_i$ for $i=1,2$.
\end{lemma}
\dem Consider a block decomposition of $R$ adapted to the one of $A$:
$$M=\left(\begin{array}{c|c} R_{1}&R_3\\ \hline R_4&R_2\end{array}\right).$$
The identity  $\varepsilon R=A^TRA$ then reads
$$\left(\begin{array}{c|c} \varepsilon R_{1}&\varepsilon R_3\\ \hline \varepsilon R_4&\varepsilon R_2\end{array}\right)=
\left(\begin{array}{c|c} A_1^TR_{1}A_1&A_1^TR_3A_2\\ \hline \stackrel{\,}{A_2^TR_4A_1}&A_2^TR_2A_2\end{array}\right).$$
Since $\mbox{\rm Sp}(A_1)\cap\mbox{\rm Sp}(\varepsilon A_2^{-1})=\emptyset$, the matrices $A_2^T\otimes A_1^T$ and $A_1^T\otimes A_2^T$ do not admit the eigenvalue $\varepsilon$ and, therefore, $R_3=0$ and $R_4=0$, from where the result follows.\qed

\medskip

The proof of next result is straightforward.

\begin{lemma} Let us choose $\varepsilon=\pm 1$, consider $\lambda\in{\mathbb C}$ such that  $\lambda^2\not\in\{ 0,\varepsilon\}$ and let $A=\mbox{\rm diag} (\lambda I,\varepsilon\lambda ^{-1} I)\in{\C}^{2m\times 2m}$. A symmetric matrix $R=R^T\in  {\C}^{2m\times 2m}$ verifies $A^TRA=\varepsilon R$ if and only if 
$$R=\left(\begin{array}{c|c} 0&B\\ \hline \stackrel{\,}{B^T}&0\end{array}\right)$$
for some $B\in {\C}^{m\times m}$.
\end{lemma}

\begin{remark}{\em The case $\lambda^2=\varepsilon$, not considered in the lemma, is also straightforward since one obviously has $\varepsilon\lambda ^{-1}=\lambda$ and, thus, $A=\mbox{\rm diag} (\lambda I,\varepsilon\lambda ^{-1} I)=\lambda I$, which clearly shows that every  matrix $R=R^T\in  {\C}^{2m\times 2m}$ verifies   $A^TRA=\varepsilon R$. It is clear that this is also the case whenever $A=\lambda I$ for $\lambda^2=\varepsilon$ and $A$ has odd order.}
\end{remark}
 
 \bigskip

\begin{theorem}\label{existence} Let us consider $\varepsilon\in\{1,-1\}$ and let $A\in\R^{(n+1)\times (n+1)}$ be a semisimple invertible matrix which is similar to $\varepsilon A^{-1}$ and define $\mathcal{C}_\varepsilon (A)$ as in Proposition \ref{red-co}. 
\begin{enumerate}
\item[(1)] There exists a symmetric invertible matrix $M\in\R^{(n+1)\times (n+1)}$ such that $A^TMA=\varepsilon M$.
\item[(2)] Let 
us consider the set
$\sigma_\varepsilon(A)$ of  all distinct eigenvalues of $A$ verifying either $|\lambda|<1$ or $|\lambda|=1$ with $(\lambda-\varepsilon\overline\lambda)>0$. If we denote the algebraic multiplicity of  each $\lambda\in\mbox{\rm Sp}(A)$ by $m(\lambda)$, then
we have
 $$\mbox{\rm dim}\, \mathcal{C}_\varepsilon(A)=\frac{m(\sqrt{\varepsilon})^2+m(\sqrt{\varepsilon})+m(-\sqrt{\varepsilon})^2+m(-\sqrt{\varepsilon})}{2}+\sum_{\lambda\in\sigma_\varepsilon(A)}m(\lambda)^2.$$
 \end{enumerate}
\end{theorem}
\dem If we consider $\sigma_\varepsilon(A)=\{\lambda_1,\dots,\lambda_k\}$ then, there exists $P\in\C^{(n+1)\times (n+1)}$ such that $A=PDP^{-1}$ for 
\begin{eqnarray*}
& & D=\mbox{diag}(D_1,\dots,D_k,\sqrt{\varepsilon} I_r,-\sqrt{\varepsilon}I_s)\, ,
\\
& & D_i=\mbox{diag}(\lambda_iI_{m(\lambda_i)},{\varepsilon}\lambda_i^{-1}I_{m(\lambda_i)})\, ,\qquad \lambda_i\in\sigma(A)\, ,
\end{eqnarray*}
where $r=m(\sqrt{\varepsilon})$ and $s=m(-\sqrt{\varepsilon})$ (which could be 0, and are equal if $\varepsilon=-1$). The congruence map $\psi (N)=P^TNP$  provides a linear isomorphim between $\mathcal{C}(A)^\C=\{N\in \C^{(n+1)\times (n+1)}\, :\, N=N^T,\ A^TNA=\varepsilon N\}$ and the complex space $\mathcal{V}=\{ R\in \C^{(n+1)\times (n+1)}\, :\, R=R^T,\ DRD=R\}.$ 

The  lemmas and the remark above show that $R\in \mathcal{V}$ if and only if $R=\mbox{diag} (R_1,\dots,R_k,S_1,S_2)$ where $S_1,S_2$ are arbitrary complex symmetric matrices of respective dimension $r,s$ and each $R_i$ is given by
$$R_i=\left(\begin{array}{c|c} 0&B_i\\ \hline \stackrel{\,}{B_i^T}&0\end{array}\right)$$
for some $B_i\in {\C}^{m(\lambda_i)\times m(\lambda_i)}$.
We then have that
$$\mbox{\rm dim}_{\C}(\mathcal{V})=m(\lambda_1)^2+\dots+m(\lambda_k)^2+\frac{r^2+r}{2}+\frac{s^2+s}{2}.$$
Further, we can choose all the matrices $S_1,S_2$ and $B_i$, $1\le i\le k$ invertible and we then get that $R$ itself is invertible. Now, both assertions of the theorem follow at once from Proposition \ref{red-co}.\qed 

\section{Orbits contained in non-degenerate invariant quadrics}
  
  Notice that in the proof of Theorem \ref{existence}, we could have considered a complex diagonal form similar to $A$ as follows:
 $$D=\mbox{\rm diag}(D(\lambda_1),\dots, D(\lambda_j), \sqrt{\varepsilon} I_r,-\sqrt{\varepsilon}I_s)$$
where $\lambda_1,\dots,\lambda_j$ are the eigenvalues in $\sigma_\varepsilon(A)$ but repeated as many times as their algebraic multiplicity, and $D(\lambda)=\mbox{diag}(\lambda,\varepsilon\lambda^{-1})$. In such case, if we consider $$K=\left(\begin{array}{cc} 0&1\\1 &0\end{array}\right),$$ then the symmetric matrix
\begin{equation}\label{matS} S=\mbox{\rm diag}(\alpha_1 K,\dots,\alpha_jK,\alpha_{j+1},,\dots ,\alpha_{n+1})\end{equation}
clealy verifies $DSD=\varepsilon S$ and, when $\alpha_i\ne 0$ for all $i\le n+1$, will provide, as in the proof of the theorem, a non-degenerate real and symmetric solution to $A^TMA=\varepsilon M$.

We will now use such matrix $S$ to prove that, under the assumptions of Theorem \ref{existence}, every non-fixed orbit is contained in an invariant quadric or an invariant affine variety.

\begin{theorem} Let $F=q\circ A\circ\ell$ where $A\in\R^{(n+1)\times(n+1)}$ is a semisimple invertible matrix which is similar to $\varepsilon A^{-1}$ for some $\varepsilon\in\{1,-1\}$.

If $X_0\not\in\mathcal{PF}$ is not contained in an $F$-invariant affine variety, then there exists at least one non-degenerate invariant quadric which contains the whole orbit of $X_0$.

\end{theorem}
\dem 
Consider a vector $X_0\not\in\mathcal{PF}$. Let us consider the complex diagonal  matrix $D$ similar to $A$ as in  the remark above and  $P$ a complex matrix such that $D=PAP^{-1}$. The vector
$$\ell(X_0)^TP^T=(a_1,b_1,\dots,a_k,b_k,a_{k+1},\dots,a_{k+r+s})\in\C^{n+1}$$
represents the coordinates of $\ell(X_0)$ with respect to the corresponding  basis of complex eigenvectors. Put $c_i=2a_ib_i$ for $i\le k$ and $c_i=a_i^2$ for $i> k$. 

If there exist $i_1\ne i_2$ such that $c_{i_1}\ne 0\ne c_{i_2}$ then we can obviously choose   complex numbers $\alpha_i\ne 0$ for all $i\ne i_1$ such that
\begin{equation}\label{inQ} \alpha_{i_1}=-c_{i_1}^{-1}\sum_{i\ne i_1}\alpha_ic_i\ne 0.
\end{equation}
Now, let us define $R=P^TSP$ where $S$ is the diagonal matrix given by (\ref{matS}) for the chosen values of $\alpha_i$. It is then clear that
$$\ell(X_0)^TR\ell(X_0)=\ell(X_0)^TP^TSP\ell(X_0)=\sum_{i=1}^{k+r+s}\alpha_ic_i=0.$$
Moreover, the invertible symmetric matrix $R$ also verifies $A^TRA=\varepsilon R$ and, hence, if $M_1,M_2$ are respectively the real and imaginary parts of $R$, then there is $\mu\in\R$ such that $M=M_1+\mu M_2$ is invertible, symmetric and  $A^TMA=\varepsilon M$. A trivial calculation shows that, also,  
$\ell(X_0)^TM\ell(X_0)=0$ because $\ell(X_0)$ is real. This means that $X_0$ lies in the $F$-invariant quadric ${\mathcal Q}(M)$ and, accordingly, so does $F^m(X_0)$ whenever it is defined.

Now suppose that we can only find a unique $i_0$ for which $c_{i_0}\ne 0$. Since $n+1=2k+r+s\ge 3$, we can find $l$ such that $c_l=0$. This implies that at least one of the coordinates of $\ell (X_0)$ with respect to the complex  basis of eigenvectors is null. Denote by $v_l$ the eigenvector corresponding to such null coordinate. If  $v_l$ is real, then  $\ell (X_0)$ is contained in the $n$-dimensional complex $A$-invariant subspace of  $\R^{n+1}$ which is supplementary to $\R$-span$\{v_l\}$, whereas if $v_l$ is not real, then the coordinate of $\ell (X_0)$ corresponding to the complex conjugate of $v_l$ must also vanish. This shows that $\ell (X_0)$ is contained in an $(n-1)$-dimensional $A$-invariant subspace of $\R^{n+1}$. In both cases we have, according to Proposition  \ref{invpl}, that $X_0$ lies on an $F$-invariant affine variety.
 \qed
 
 \bigskip
  
\begin{remark} {\em  It is interesting to point out the following:
\begin{enumerate}
\item[(1)] When for a given $X_0\in\R^n$, we have that equation (\ref{inQ})  holds for different choices of $\alpha_i$, $i\ne i_1$, the point $X_0$ is contained in the corresponding quadrics. This is always the case if one may find three distinct $i_1,i_2,i_3$ such that $c_{i_l}\ne 0$ for $l=1,2,3$. This means that, generally speaking, orbits are usually contained in the intersection of two or more non-degenerate quadrics.
\item[(2)] The case in which all the eigenvalues verify $\lambda^2=\epsilon$ is  globally $2$-periodic \cite{b-l}, this meaning that all the orbits are either fixed points or of the form $\{X_0,X_1,X_0,X_1\cdots\}$. It is then obvious that every non fixed orbit lies on the unique line connecting $X_0$ and $X_1$. Further,  such line is always $F$-invariant. To see this it suffices to show that $U=\R\mbox{-span}\{\ell(X_0),\ell(X_1)\}$ is $A$-invariant, which is obvious since for $i,j\in\{0,1\}$, $i\ne j$ we have
$q(\ell(X_i))=X_i=F(X_j)=q(A\ell(X_j))$ and this implies that there exists $\alpha\ne 0$ such that $\ell(X_i)=A\ell(X_j)$.
\end{enumerate}}
\end{remark}

We finish this section with a quite relevant result in the cases in which $A$ is not similar to an orthogonal matrix. 
It is not difficult to see that if the matrix $A$ has an equilibrium associated to an eigenvalue of maximal modulus with algebraic multiplicity 1, then it atracts almost all solutions. This means that the such fixed point should be adherent to almost all invariant quadrics. Actually, we have the following stronger result:

\begin{proposition} \label{inc} Let $F=q\circ A\circ\ell$ where $A\in\R^{(n+1)\times(n+1)}$ is a semisimple invertible matrix which is similar to $\varepsilon A^{-1}$ for some $\varepsilon\in\{1,-1\}$ and suppose that $\lambda\in\mbox{\rm Sp}(A)$ verifies $|\lambda|\ne 1$.
\begin{enumerate}
\item[(1)] If $\lambda\in\R$ and $X_\lambda$ is a fixed point of $F$ associated  to $\lambda$, then
every $F$-invariant non-degenerate quadric contains the point $X_\lambda$.
\item[(2)] If $\lambda\not\in\R$ and $v\in\C^{n+1}$ is an associate complex eigenvector such that the real linear span $U$ of its real and imaginary parts is not contained in $U_0$, then every $F$-invariant non-degenerate quadric contains the $F$-invariant line in $\R^n$ defined by $U$.
\end{enumerate}
\end{proposition}
\dem Consider an $F$-invariant non-degenerate quadric ${\mathcal Q}(M)$. Since it holds $A^TMA=\varepsilon M$, where $\varepsilon =\pm 1$, for each complex vector $v\in\C^{n+1}$ such that $Av=\lambda v$ we get that $|\lambda|^2v^TM\overline v=\varepsilon v^TM\overline v$ and, therefore, $v^TM\overline v=0$ because $|\lambda|\ne 1$. If $\lambda\in\R$ and $X_\lambda\in\R^n$ is a fixed point of $F$, then the result follows at once by considering $v=\ell (X_\lambda)$. Suppose now that $\lambda\not\in\R$ and $X_0\in\R^n$ lies in the $F$-invariant line defined by $U=\R\mbox{-span}\{v+\overline v,i(v-\overline v)\}$. This means that $\ell (X_0)=\alpha v+\beta \overline v$ for appropriate values of $\alpha,\beta\in\C$. It is clear that from $A^TMA=\varepsilon M$ one also gets
$\lambda^2v^TM v=\varepsilon v^TM v$ and $\overline\lambda^2\overline v^TM \overline v=\varepsilon\overline v^TM \overline v$. Since $\lambda^2\ne \varepsilon\ne \overline\lambda ^2$, this implies $v^TM v=\overline v^TM \overline v=0$ and the result follows at once.\qed

 \begin{remark} {\em  The same proof shows that the results of Proposition \ref{inc} are still valid in the case $\varepsilon=-1$ and $\lambda^2\ne -1$.
}\end{remark}

\section*{Acknowledgements} 

The author is in debt with Prof. Eduardo Liz for his useful suggestions and the careful reading of the paper.
This work was partially supported by  the Spanish Ministry of Science and Innovations and FEDER, grant MTM2010-14837.

\

\noindent {\it Author's address:}\\
Depto. Matem\'atica Aplicada II, E.I.Telecomunicaci\'on, Universidad de Vigo, 36310 Vigo, Spain\\
{ibajo@dma.uvigo.es}

\end{document}